\pdfoutput=1
\documentclass[a4paper,11pt]{amsart}
\usepackage{amssymb}
\usepackage{amsfonts}
\usepackage{amsmath}
\usepackage{mathpazo}
\usepackage{hyperref}
\usepackage{pdfpages}
\usepackage{tabularx}
\usepackage{colortbl}
\usepackage{multimedia}
\usepackage{rotating}
\usepackage{tikz}
\usepackage{epstopdf}
\usepackage{graphicx}
\usepackage{rotating}
\usepackage{enumitem}
\usepackage{color}
\usepackage{listings}
\usepackage[noend]{algorithmic}%

\lstdefinelanguage{Julia}%
  {morekeywords={abstract,break,case,catch,const,continue,do,else,elseif,%
      end,export,false,for,function,immutable,import,importall,if,in,%
      macro,module,otherwise,quote,return,switch,true,try,type,typealias,%
      using,while,Integer},%
   sensitive=true,%
   alsoother={$},%
   morecomment=[l]\#,%
   morecomment=[n]{\#=}{=\#},%
   morestring=[s]{"}{"},%
   morestring=[m]{'}{'},%
}[keywords,comments,strings]%

\usepackage{graphicx}
\usepackage{url}

\newtheorem{theorem}{Theorem}

\newtheorem{definition}[theorem]{Definition}
\newtheorem{lemma}[theorem]{Lemma}
\newtheorem{algorithm}[theorem]{Algorithm}
\theoremstyle{remark}
\newtheorem{example}[theorem]{Example}

\begin{document}

\title{Bad Primes in Computational Algebraic Geometry}  

\author[J.~B\"ohm]{Janko~B\"ohm}
\address{Department of Mathematics\\
University of Kaiserslautern\\
Erwin-Schr\"odinger-Str.\\
67663 Kaiserslautern\\
Germany}
\email{boehm@mathematik.uni-kl.de}

\author[W.~Decker]{Wolfram Decker}
\address{Department of Mathematics\\
University of Kaiserslautern\\
Erwin-Schr\"odinger-Str.\\
67663 Kaiserslautern\\
Germany}
\email{decker@mathematik.uni-kl.de}

\author[C.~Fieker]{Claus Fieker}
\address{Department of Mathematics\\
University of Kaiserslautern\\
Erwin-Schr\"odinger-Str.\\
67663 Kaiserslautern\\
Germany}
\email{fieker@mathematik.uni-kl.de}

\author[S. Laplagne]{Santiago Laplagne}
\address{Santiago Laplagne\\
Departamento de Matem\'atica\\
FCEN\\
Universidad de Buenos Aires - Ciudad Universitaria\\
Pabell\'on I - (C1428EGA) - Buenos Aires\\
Argentina}
\email{slaplagn@dm.uba.ar}

\author[G.~Pfister]{Gerhard Pfister}
\address{Department of Mathematics\\
University of Kaiserslautern\\
Erwin-Schr\"odinger-Str.\\
67663 Kaiserslautern\\
Germany}
\email{pfister@mathematik.uni-kl.de}

\keywords{Modular computations, bad primes, algebraic curves, adjoint ideal}

\maketitle
\begin{abstract}
Computations over the rational numbers often suffer from intermediate
coefficient swell.  
 One solution to this problem is to apply the given algorithm
modulo a number of primes and then lift the modular results to the rationals. This method is
guaranteed to work if we use a sufficiently large set of good primes. In many
applications, however, there is no efficient way of excluding bad primes. In this note, we
describe a technique for rational reconstruction which will nevertheless
return the correct result, provided the number of good primes in the selected set of primes is
large enough.   We give a number of illustrating examples which are implemented using the computer algebra system \textsc{Singular} and the programming language \textsc{Julia}. We discuss applications of our technique in computational
algebraic geometry.
\end{abstract}

\section{Introduction}
Many exact computations in computer algebra are carried out over the rationals
and extensions thereof. Modular techniques are an important tool to improve
the performance of such algorithms since intermediate coefficient growth is
avoided and the resulting modular computations can be done in parallel. For this, we require that the
algorithm under consideration is also applicable over finite fields and returns a
deterministic result. The fundamental approach is then as follows: Compute the
result modulo a number of primes. Then reconstruct the result over $\mathbb{Q}$
from the modular results. 
\begin{example}\label{ex1}
To compute%
\[
\frac{1}{2}+\frac{1}{3}=\frac{5}{6}%
\]
using modular methods, the first step is to apply the Chinese remainder isomorphism:

\[%
\begin{tabular}
[c]{ccccccccccc}
&  & $\mathbb{Z}/5$ & $\times$ & $\mathbb{Z}/7$  & $\times$ & $\mathbb{Z}/101$ & $\cong$ & $\mathbb{Z}/3535\medskip$\\
$\frac{1}{2}$ & $\longmapsto$ & $($ $\overline{3}$ & $,$ & $\overline{4}$ & $,$ & $\overline{51}$ $)$ &  & $\smallskip$\\
&  &  &  &$+$ &  &  &  &  & $\smallskip$\\
$\frac{1}{3}$ & $\longmapsto$ & $($ $\overline{2}$ & $,$ & $\overline{5}$ & $,$ & $\overline{34}$ $)$ &  & $\smallskip$\\
&  &  &  &$\shortparallel$ &  &  &  &  & $\smallskip$\\
&  & $($ $\overline{0}$ & $,$ & $\overline{2}$ & $,$ & $\overline{85}$ $)$ & $\longmapsto$ & $\overline{590}$%
\end{tabular}
\]

\end{example}
The second step is to reconstruct a rational number from $\overline{590}$.

\section{Rational Reconstruction}

\begin{theorem}\cite{KG}
For every integer $N$, the $N$-\textbf{Farey map}%
\[
\hspace{-3mm}%
\begin{tabular}
[c]{rcc}%
$\left\{  \frac{a}{b}\in\mathbb{Q}\left\vert
\begin{tabular}
[c]{r}%
$\gcd(a,b)=1$\\
$\gcd(b,N)=1$%
\end{tabular}
\ \ \right.  \text{ }\left\vert a\right\vert ,\left\vert b\right\vert
\leq\sqrt{(N-1)/2}\right\}  $ & $\longrightarrow$ & $\mathbb{Z}/N\medskip$\\
$\frac{a}{b}$ & $\longmapsto$ & $\overline{a}\cdot\overline{b}^{-1}$%
\end{tabular}
\ \ \
\]
is injective.
\end{theorem}

There are efficient algorithms for computing preimages of the Farey map, see, for example, \cite[Sec. 5]{KG}.

\begin{example}
We use the computer algebra system \textsc{Singular} \cite{DGPS} to compute the preimage of the Farey map the setting of Example \ref{ex1}:
\smallskip

\noindent { \texttt{$\color{blue}>$ \color{red}ring \color{black}r = 0, x, dp;}}

\noindent { \texttt{$\color{blue}>$ \color{red}farey\color{black}(590,3535);}}

\noindent { \texttt{\color{blue} \hspace{0.3cm} 5/6}}

\end{example}

The basic concept for modular computations is then as follows:

\begin{enumerate}
\item Compute the result over $\mathbb{Z}/p_{i}$ for distinct primes $p_{1}%
,\ldots,p_{r}$.

\item 
Use the Chinese remainder isomorphism
\[
\mathbb{Z}/N\hspace{2mm}\cong\hspace{2mm}\mathbb{Z}/p_{1}\times\ldots
\times\mathbb{Z}/p_{r}%
\]
to lift the modular results to  $\mathbb{Z}/N$ where $N=p_{1}\cdots p_{r}$.

\item Compute the preimage of the lift with respect to the $N$-Farey map.

\item Verify the correctness of the lift.
\end{enumerate}

This will yield the correct result, provided $N$ is large enough (that is, the
$\mathbb{Q}$-result is contained in the domain of the $N$-Farey map), and provided none of the $p_{i}$ is bad.

\begin{definition}
A prime $p$ is called \textbf{bad} (with respect to a fixed algorithm and input) if the
result over $\mathbb{Q}$ does not reduce modulo $p$ to the result over
$\mathbb{Z}/p$.
\end{definition}
By convention, this includes the case where, modulo $p$,  the input is not defined or the algorithm in consideration is not applicable.

\section{Bad primes}

\subsection{Bad primes in Gr\"{o}bner basis computations}

Consider a set of variables $X=\{x_1,\ldots,x_n\}$ and a monomial ordering $>$ on
the monomials in $X$. For a set of polynomials $G$, write
$\operatorname{LM}(G)$ for its set of lead monomials. For $G\subset
\mathbb{Z}[X]$ and $p$ prime, write $G_{p}$ for the image of $G$ in $\mathbb{Z}\slash\hspace
{-0.5mm} p\hspace
{0.25mm}[X]$.

\begin{theorem}\cite{arnold}
 Suppose $F=\{f_{1},...,f_{r}\}\subset\mathbb{Z}[X]$ with all $f_{i}$
primitive and homogeneous. Let $G$ be the reduced Gr\"{o}bner basis of $\left\langle F\right\rangle
\subset\mathbb{Q}[X]$,
$G(p)$ the reduced Gr\"{o}bner basis of $\left\langle F_{p}%
\right\rangle $, and
$G_{\mathbb{Z}}$ a minimal strong Gr\"{o}bner basis of $\left\langle
F\right\rangle \subset\mathbb{Z}[X]$. Then 
\begin{center}
$p$ does not divide any lead coefficient in $G_{\mathbb{Z}}$
$\Leftrightarrow$ $\operatorname{LM}(G)=\operatorname{LM}(G(p))$ $\Leftrightarrow$ $ G_{p}=G(p)$.
\end{center}
\end{theorem}

\begin{example} Using \textsc{Singular}, we determine the bad primes for a Gr\"obner basis computation of the Jacobian ideal of a projective plane curve. We compute a minimial strong Gr\"obner basis over $\mathbb{Z}$:\smallskip

\noindent {\texttt{$\color{blue}>$ \color{red} option\color{black}("redSB");}}\\
\noindent { \texttt{$\color{blue}>$ \color{red} ring \color{black}R = \color{black}integer\color{black},(x, y, z),lp;}}\\
\noindent { \texttt{$\color{blue}>$ \color{red} poly \color{black}f = x7y5 + x2yz9 + xz11 + y3z9;}}\\
\noindent { \texttt{$\color{blue}>$ \color{red} ideal \color{black}I = \color{red}groebner\color{black}(\color{red}ideal\color{black}(\color{red}diff\color{black}(f, x), \color{red}diff\color{black}(f, y), \color{red}diff\color{black}(f,z)));}}\\
\noindent { \texttt{$\color{blue}>$  \color{red} apply(\color{red}list\color{black}(I[1..\color{red}size\color{black}(I)]),\color{red}leadcoef\color{black})\color{black};}}\smallskip\\
\noindent { \texttt{\color{blue} 13781115527868730344777310464613260 83521912290113517241074608876444 60 12 4 12 12 45349632 12 1473863040 12 22674816 12 3888 12 12 12 13608 12 108 54 6 2 27 3 1 4 2 2 1 216 1 2 3 1 540 12 108 27 3 1 9 3 1 1 1 1 1 7 1 5 1 1}}

\bigskip
\noindent The bad primes, that is, the primes $p$ with $G_{p}\neq G(p)$, are then the prime factors
$$ p=2,3,5,7,11,13,257,247072949,328838088993550682027
$$
of the lead coefficients. In contrast, the lead coefficients of the Gr\"obner basis
over $\mathbb{Q}$ involve only the prime factors $2,3,5,7,13$, and hence not all bad primes. As shown by the following computation, $257$ is indeed a bad prime:
\bigskip

\noindent { \texttt{$\color{blue}>$ \color{red} ring \color{black}R0 = 0,(x, y, z),lp;}}\\
\noindent { \texttt{$\color{blue}>$ \color{red}size\color{black}(\color{red}lead\color{black}(\color{red}groebner\color{black}(\color{red}fetch\color{black}(R,I))));}}\\
\noindent { \texttt{$\color{blue} 15$}}\\
\noindent { \texttt{$\color{blue}>$ \color{red} ring \color{black}R1 = 257,(x, y, z),lp;}}\\
\noindent { \texttt{$\color{blue}>$ \color{red}size\color{black}(\color{red}lead\color{black}(\color{red}groebner\color{black}(\color{red}fetch\color{black}(R,I))));}}\\
\noindent { \texttt{$\color{blue} 14$}}
\end{example}

\subsection{Classification of Bad Primes}

Bad primes can be classified as follows, see \cite[Sec. 3]{BDFP} for details:
\begin{itemize}[leftmargin=*]
\item Type 1: The input modulo $p$ is not valid (this poses no problem).
\item Type 2: There is a failure in the course of the algorithm (for example, a matrix may not be
invertible modulo $p$; this wastes computation time if it happens).
\item Type 3: A computable invariant with known expected value (for example, a Hilbert
polynomial) has a wrong value in a modular computation (to detect this we have to do expensive tests for each prime, although the set
of bad primes usually is finite, and hence bad primes rarely occur).
\item Type 4: A computable invariant with unknown expected value (for example, the lead
ideal in a Gr\"{o}bner basis computation) is wrong (this can be handled by a majority
vote, however we have to compute the invariant for each modular result and store the modular results).
\item Type 5: otherwise.
\end{itemize}

The Type 5 case in fact occurs, as is shown by the following example. For an ideal $I\subset\mathbb{Q}[X]$ and a prime $p$ define $I_{p}=(I\cap
\mathbb{Z}[X])_{p}$.
\begin{example}
Consider the algorithm $I\mapsto\sqrt{I+\operatorname{Jac}(I)}$ computing the radical of the Jacobian ideal for the curve%
\[
\ I=\scalebox{1}{$\left\langle x^{6}+y^{6}+7x^{5}z+x^{3}y^{2}z-31x^{4}z^{2}-224x^{3}z^{3}+244x^{2}z^{4}+1632xz^{5}+576z^{6}\right\rangle$.}%
\]
Note that, with respect to the degree reverse lexicographic order, $\operatorname{LM}(I)=\left\langle x^{6}\right\rangle =\operatorname{LM}%
(I_{5})$, that is, $5$ is not bad with respect to the input. The following computation in \textsc{Singular} first determines the minimal associated primes of $U(0)  =\sqrt{I+\operatorname{Jac}(I)}$ and $U(5)  =\sqrt{I_5+\operatorname{Jac}(I_5)}$.\smallskip

\noindent {\texttt{$\color{blue}>$  \color{red} LIB \color{black} "primdec.lib";}}\\
\noindent { \texttt{$\color{blue}>$  \color{red} ring \color{black} R0 = 0, (x, y, z), dp;}}\\
\noindent { \texttt{$\color{blue}>$  \color{red} poly \color{black} f = x6+y6+7x5z+x3y2z-31x4z2-224x3z3+244x2z4+1632xz5+576z6;}}\\
\noindent { \texttt{$\color{blue}>$  \color{red}  ideal \color{black} U0 = \color{red}radical\color{black}(\color{red}ideal\color{black}(f, \color{red}diff\color{black}(f, x), \color{red}diff\color{black}(f, y), \color{red}diff\color{black}(f, z)));}}\\
\noindent { \texttt{$\color{blue}>$  \color{red} minAssGTZ\color{black}(U0);}}\\
\noindent { \texttt{\color{blue}[1]:  \underline{\phantom{w}}[1]=y\hspace{1cm}[2]:  \underline{\phantom{w}}[1]=y\\}}
\noindent {\texttt{\phantom{[1]:\,\,}\color{blue}\underline{\phantom{w}}[2]=x+6z \hspace{1.175cm}\underline{\phantom{w}}[2]=x-4z}}\\
\noindent { \texttt{$\color{blue}>$  \color{red} ring \color{black} R5 = 5, (x, y, z), dp;}}\\
\noindent { \texttt{$\color{blue}>$  \color{red} poly \color{black} f =\color{red}imap\color{black}(R0,f);}}\\
\noindent { \texttt{$\color{blue}>$  \color{red}  ideal \color{black} U5 = \color{red}radical\color{black}(\color{red}ideal\color{black}(f, \color{red}diff\color{black}(f, x), \color{red}diff\color{black}(f, y), \color{red}diff\color{black}(f, z)));}}\\
\noindent { \texttt{$\color{blue}>$  \color{red} minAssGTZ\color{black}(U5);}}\\
\noindent { \texttt{\color{blue}[1]:  \underline{\phantom{w}}[1]=y\hspace{1cm}[2]:  \underline{\phantom{w}}[1]=y\\}}
\noindent {\texttt{\phantom{[1]:\,\,}\color{blue}\underline{\phantom{w}}[2]=x-z \hspace{1.3cm}\underline{\phantom{w}}[2]=x+z}}\\
\noindent { \texttt{$\color{blue}>$  \color{red}minassGTZ\color{black}(\color{red}imap\color{black}(R0,U0));}}\\
\noindent { \texttt{\color{blue}[1]: \underline{\phantom{w}}[1]=y}}\\
\noindent { \texttt{\color{blue}\phantom{[1]: }\underline{\phantom{w}}[2]=x+z}}

\bigskip
\noindent This shows that $U(0)_{5}\neq U(5)$, but $\operatorname{LM}(U(0))=\left\langle y,x^{2}\right\rangle =\operatorname{LM}%
(U(5))$.
\end{example}

\section{Error-Tolerant Reconstruction}

Our goal is to reconstruct the $\mathbb{Q}$-result $\frac{a}{b}$ from the modular result $\overline{r}\in\mathbb{Z}/N$ in
the presence of bad primes. Our basic strategy will be to find an element $(x,y)$ with $\frac{x}{y}=\frac{a}{b}$ in the
lattice%
\[
\Lambda=\langle(N,0),(r,1)\rangle\subset\mathbb{Z}^{2}.%
\]

\begin{lemma}\label{lem collinear}
\cite[Lem. 4.2]{BDFP} All $(x,y)\in\Lambda$ with $x^{2}+y^{2}<N$ are collinear.
\end{lemma}

Now suppose $N=N^{\prime}\cdot M$
with $\gcd(N^{\prime},M)=1$. We assume that $N^{\prime}$ is the product of the good
primes with correct result $\overline{s}$, and $M$ is the product of the
bad primes with wrong result~$\overline{t}$.

\begin{theorem}
\cite[Lem. 4.3]{BDFP} If%
\[
\overline{r}\mapsto(\overline{s},\overline{t})\text{ \quad with respect
to}\quad\mathbb{Z}/N\cong\mathbb{Z}/N^{\prime}\times\mathbb{Z}/M
\]
and%
\[
\frac{a}{b} \operatorname{mod}N^{\prime} = s%
\]
then $(aM,bM)\in\Lambda$. So if $(a^{2}+b^{2})M<N^{\prime}$,
then (by Lemma \ref{lem collinear})%
\[
\frac{x}{y}=\frac{a}{b}\hspace{2mm}\text{ for all }(x,y)\in\Lambda\text{ with
}(x^{2}+y^{2})<N%
\]
and such vectors exist. Moreover, if $\gcd(a,b)=1$ and $(x,y)$ is a shortest
vector $\neq0$ in $\Lambda$, we also have $\gcd(x,y)|M$.
\end{theorem}

Hence, if $N^{\prime}\gg M$, the Gauss-Lagrange-Algorithm for finding a
shortest vector $(x,y)\in\Lambda$ gives $\frac{a}{b}$ independently of $t$,
provided $x^{2}+y^{2}<N$. We use the programming language \textsc{Julia}\footnote{See \url{http://julialang.org/}.},
to illustrate the resulting algorithm. 

\begin{lstlisting}
function ErrorTolerantReconstruction(r::Integer, N::Integer)
  a1 = [N, 0]
  a2 = [r, 1]
  while dot(a1, a1) > dot(a2, a2)
    q = dot(a1, a2)//dot(a2, a2)
    a1, a2 = a2, a1 - Integer(round(q))*a2
  end
  if dot(a1, a1) < N
    return a1[1]//a1[2]
  else
    return false
  end
end
\end{lstlisting}

The following table shows timings (in seconds), for $r$ and $N$ of bit-length $500$, comparing the \textsc{Julia}-function with implementations in
the \textsc{Singular}-kernel (optimized C/C++ code) and the current \textsc{Singular}-interpreter:
\begin{center}
\def\arraystretch{1.1}
\begin{tabular}{cc|ccc|cc}
\textsc{Singular}-kernel &&&\textsc{Julia} &&& \textsc{Singular}-interpreter \\\hline
0.001&&&0.005 &&& 0.055  
\end{tabular}
\end{center}
Building on \textsc{Julia} as a fast mid-level language, a backwards-compatible just-in-time compiled \textsc{Singular}-interpreter is under development. 

\begin{example}
In the setting of Example \ref{ex1}, we obtain $\frac{5}{6}$ from
$\overline{590}\in\mathbb{Z}/3535$ by%
\begin{lstlisting}
julia> ErrorTolerantReconstruction(590, 3535)
5//6
\end{lstlisting}
which computes the sequence
\begin{align*}
(3535,0)  &  =6\cdot(590,1)+(-5,-6),\\
(590,1)  &  =-48\cdot(-5,-6)+(350,-287).
\end{align*}

\end{example}

\begin{example}
Now we introduce an \color{red}error \color{black}in the modular results:
\[%
\begin{tabular}
[c]{ccccccccc}%
$\mathbb{Z}/5$ & $\times$ & $\mathbb{Z}/7$  &
$\times$ & $\mathbb{Z}/101$ & $\cong$ & $\mathbb{Z}/3535\medskip$\\
$($ \color{red}$\overline{1}$\color{black} & $,$ &
$\overline{2}$ &  & $\overline{85}$ $)$ & $\mapsto$ & $\overline{2711}$%
\end{tabular}
\ \
\]
Error tolerant reconstruction computes%
\begin{align*}
(3535,0)  &  = 1\cdot(2711,1)+(824,-1),\\
(2711,1)  &  =3\cdot(824,-1)+(239,4)\\
(824,-1)  &  =  3  \cdot (239,4) +(107,-13)\\
(239,4) &=  2  \cdot (107,-13) + (25,30)\\
(107,-13)& =  1  \cdot (25,30) +(82,-43)
\end{align*}
hence yields%
\[
\frac{25}{30}=\frac{5\cdot5}{5\cdot6}=\frac{5}{6}\text{.}
\]
Note that\vspace{-1mm}
\[
(5^{2}+6^{2})\cdot5 = 305 < 707=7\cdot101.%
\]

\end{example}

\section{General Reconstruction Scheme for Commutative Algebra}

 For a given ideal $I\subset\mathbb{Q}[X]$, we want to compute some ideal (or module)
$U(0)$ associated to $I$ by a deterministic algorithm. We proceed along the following lines:%
\begin{enumerate}
\item Over $\mathbb{Z}/p$ compute  $U(p)$ from $I_{p}$ for $p$ in a suitable finite
set $\mathcal{P}$ of primes.

\item Replace $\mathcal{P}$ by a subset according to a majority vote on $\operatorname{LM}%
(U(p))$ (see also \cite[Rmk. 5.7]{BDFP}).

\item For $N=\prod\nolimits_{p\in\mathcal{P}}p$ compute the coefficient-wise CRT--lift
$U(N)$ to $\mathbb{Z}/N$, identifying generators by their lead monomials.

\item Lift $U(N)$ by error tolerant rational reconstruction to $U$. 

\item Test $U_{p}=U(p)$ for some random extra prime $p$.

\item Verify $U=U(0)$.

\item If the lift, test or verification fails, then enlarge $\mathcal{P}$ and repeat.
\end{enumerate}

\begin{theorem}\cite[Lem. 5.6]{BDFP} If the bad primes form a Zariski closed proper subset of
$\operatorname{Spec}\mathbb{Z}$, then this strategy terminates with the
correct result.%
\end{theorem}

\section{Computing Adjoint Ideals}

\begin{figure}
\begin{center}
\includegraphics[
height=2.2in,
width=4in
]%
{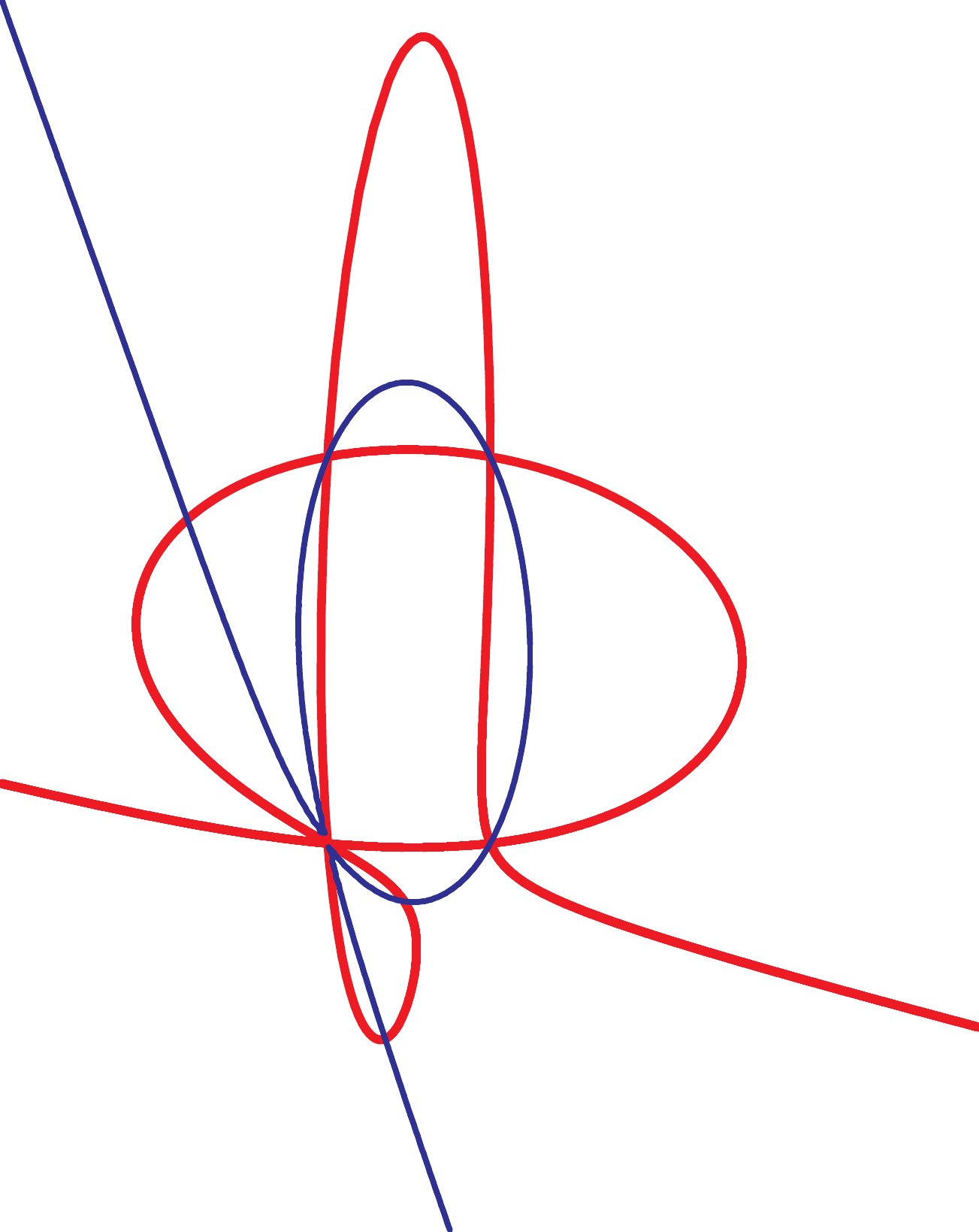}%
\end{center}
\caption{Degree $3$ adjoint curve of a rational
curve of degree $5$}
\label{adjoint}
\end{figure}

We discuss an application from algebraic geometry. The goal is to compute adjoint curves, that is, curves which pass with sufficiently high multiplicity through the singularities of a given curve, see Figure \ref{adjoint}. We consider an integral,
non-degenerate projective curve $\Gamma\subset\mathbb{P}^{r}$ with
normalization map $\pi:\overline{\Gamma}\rightarrow\Gamma$, and a saturated
homogeneous ideal $I$ with $I(\Gamma)\subsetneqq I\subset k[x_{0},...,x_{r}]$. We write $\operatorname*{Sing}(\Gamma)$ for the singular locus of $\Gamma$.
Let $H$ be the pullback of a hyperplane, and $\Delta(I)$ the pullback of
$\operatorname*{Proj}(S/I)$. Then the exact sequence%
\[
0\rightarrow\widetilde{I}\mathcal{O}_{\Gamma}\rightarrow\pi_{\ast}%
(\widetilde{I}\mathcal{O}_{\overline{\Gamma}})\rightarrow\mathcal{F}%
\rightarrow0
\]
induces, for $m\gg0$, an exact sequence%
\[
0\rightarrow I_{m}/I(\Gamma)_{m}\overset{\overline{\varrho_{m}}}{\rightarrow
}H^{0}\left(  \overline{\Gamma},\mathcal{O}_{\overline{\Gamma}}\left(
mH-\Delta(I)\right)  \right)  \rightarrow H^{0}\left(  \Gamma,\mathcal{F}%
\right)  \rightarrow0\text{.}
\]

\begin{definition}
The ideal $I$ is an \textbf{adjoint ideal} of $\Gamma$ if $\overline{\varrho_{m}}$
is surjective for $m\gg0$.
\end{definition}

Since $
h^{0}\left(  \Gamma,\mathcal{F}\right)  =\sum\nolimits_{P\in
\operatorname*{Sing}(\Gamma)}\operatorname{length}(I_{P}\overline{\mathcal{O}_{\Gamma,P}}%
/I_{P})$,
we obtain:

\begin{theorem}
\cite{AC}
With notation as above:
\[
I\text{ is an adjoint ideal of $\Gamma$ }\Longleftrightarrow\text{ }I_{P}\overline{\mathcal{O}%
_{\Gamma,P}}=I_{P}\text{ for all }P\in\operatorname*{Sing}(\Gamma).
\]
\end{theorem}
The conductor $\mathcal{C}_{\mathcal{O}_{\Gamma,P}}$ of $\mathcal{O}_{\Gamma,P}\subset\overline{\mathcal{O}%
_{\Gamma,P}}$ is the largest ideal of $\mathcal{O}_{\Gamma,P}$ which is also an ideal in $\overline{\mathcal{O}%
_{\Gamma,P}}$.

\begin{definition}The
\textbf{Gorenstein adjoint ideal} of $\Gamma$ is the largest homogeneous ideal
$\mathfrak{G}\subset K[x_{0},\ldots,x_{r}]$ with%
\[
\mathfrak{G}_{P}=\mathcal{C}_{\mathcal{O}_{\Gamma,P}}\;\text{ for all }%
\;P\in\operatorname*{Sing}(\Gamma)\text{.}%
\]

\end{definition}

The Gorenstein adjoint ideal has many applications in the geometry of curves. 

\begin{example}
If $\Gamma$ be an irreducible plane algebraic curve of degree $n$, then $\mathfrak{G}_{n-3}$ cuts out the
canonical linear series.
\end{example}
\begin{example}
If $\Gamma$ is a rational plane curve of degree $n$, then $\mathfrak{G}_{n-2}$ maps
$\Gamma$ to a rational normal curve of degree $n-2$ in $\mathbb{P}^{n-2}$.
\end{example}
\begin{example}
The Gorenstein adjoint ideal can be used in the Brill-Noether-Algorithm to
compute Riemann-Roch spaces for singular curves.
\end{example}

The Gorenstein adjoint ideal can be computed via a local-to-global strategy.

\begin{definition}
The \textbf{local adjoint ideal} of $\Gamma$ at $P\in\operatorname{Sing}%
\Gamma$ is the largest homogeneous ideal $\mathfrak{G}(P)\subset
k[x_{0},\ldots,x_{r}]$ with $\mathfrak{G}(P)_{P}=\mathcal{C}_{\mathcal{O}_{\Gamma,P}}.$

\end{definition}

\begin{lemma}
\cite[Prop. 5.4]{BDLP1}\label{lem localglobal} With notation as above,
\[
\mathfrak{G}={\bigcap\nolimits_{P\in\operatorname{Sing}\Gamma}}\mathfrak{G}(P).
\]

\end{lemma}

\begin{definition}
Let $A$ be the coordinate ring of an affine model $C%
=\operatorname{Spec}A$ of $\Gamma$ and let $P\in\operatorname{Sing}(A)$. A ring
$A\subset B\subset\overline{A}\subset \operatorname{Quot}(A)$
 is called a \textbf{minimal
local contribution} to $\overline{A}$ at $P$ if $B_{P} =\overline{A_{P}}$ and $
B_{Q}   =A_{Q}$ for all $P\neq Q\in C$.
\end{definition}

The minimal local contribution to $\overline{A}$ at $P$ is unique and can be computed using
Grauert-Remmert-type normalization algorithms, see \cite{BBDLPSS}. It can be
written as $B=\frac{U}{d}$ with an ideal $U\subset A$ and a common denominator
$d\in A$.

\begin{algorithm}
\cite[Alg. 4]{BDLP1} With notation as above, $\mathfrak{G}(P)\subset k[x_{0},\ldots,x_{r}]$ is the homogenization of the
preimage of $\left(d:U\right)$ under $ k[x_1,\ldots,x_r]\rightarrow  k[x_1,\ldots,x_r]/I=A$.
\end{algorithm}

\section{Modular version of the algorithm}

Applying the general modular strategy gives an algorithm which is  two-fold parallel (taking Lemma \ref{lem localglobal} into account). We
use primes $p$ such that the algorithm is applicable to the variety $\Gamma_{p}$ defined by
$I(\Gamma)_{p}$. Efficient verification can be realized through a semi-continuity argument, see \cite[Theorem 8.14]{BDLP1}. Table~\ref{timings} gives timings (in seconds on a 2.2 GHz processor) for plane curves $f_{n}$ of degree $n$
with $\binom{n-1}{2}$ singularities of type $A_{1}$.

 Rows
\texttt{LA} and \texttt{IQ} refer to global computations of the Gorenstein adjoint ideal via
linear algebra \cite{Mnuk} and ideal quotients, respectively. The row \texttt{Maple-IB} shows timings for the normalization of the curve
via a computation of an integral basis in \textsc{Maple} \cite{vanHoeijIntegralBasis}. The row \texttt{locIQ} gives timings
for the local-to-global (Lemma \ref{lem localglobal}), and \texttt{modLocIQ} for the modular local-to-global
strategy. We also give timings for parallel computations and for the modular probabilistic algorithm obtained by omitting the verification.
 In square brackets, the number of primes in the modular strategy is
shown, in round brackets the number of cores used simultaneously in a parallel computation.

\begin{table}[h]
\begin{center}
\def\arraystretch{1.1}
\begin{tabular}
[c]{l|c|c|rc|rc|rc|}
& {\small parallel} &
{\small probabilistic} & \text{$f_{5}$}
&  & \text{$f_{6}$} &  & \text{$f_{7}$} & \\\hline
{\small {\texttt{Maple-IB}}} &  &  & 5.1 &  & 47 &  & 318 & \\\hline
{\small {\texttt{LA}}} &  &  & 98 &  & 4400 &  & - & \\
{\small {\texttt{IQ}}} &  &  & 1.3 &  & 54 &  & 3800 & \\
{\small {\texttt{locIQ}}} & \scalebox{0.75}{$\blacksquare$} &  & 1.3 & (1) & 54 & (1) & 3800 &
(1)\\\hline
{\small {\texttt{modLocIQ}}} &  &  & 6.4 & [33] & 19 & [53] & 150 & [75]\\
&  & \scalebox{0.75}{$\blacksquare$} & 6.2 & [33] & 18 & [53] & 104 & [75]\\
& \scalebox{0.75}{$\blacksquare$} &  & .36 & (74) & 1.6 & (153) & 51 & (230)\\
& \scalebox{0.75}{$\blacksquare$} & \scalebox{0.75}{$\blacksquare$} & .21 & (74) & 0.48 & (153) & 5.2 & (230)
\end{tabular}
\end{center}

\medskip
\caption{Timings}
\label{timings}
\end{table}

Observe that, in the example, a local-to-global strategy does not give
any benefit when computing over the rationals, since the singular locus does
not decompose. However, by Chebotarev's density theorem, the singular locus
is likely to decompose when passing to a  finite field, as illustrated by the last two rows of
the table.


\vspace{3cm}
\begin{thebibliography}{44}        
\bibitem{AC}Arbarello, E.;
Ciliberto, C.: \emph{Adjoint hypersurfaces to curves in }$\mathbb{P}^{r}%
$\emph{ following Petri}, in Commutative Algebra, Lecture Notes in Pure and
Applied Mathematics, vol. 84, Dekker, New York, 1-21 (1983).
                                                                                  %
\bibitem{arnold}Arnold, E. A.: \emph{Modular algorithms for
computing Gr\"{o}bner bases}, J. Symb. Comput. 35, 403--419 (2003).


\bibitem{BDFP}J.~B\"{o}hm, W.~Decker, C.~Fieker, G.~Pfister. \emph{The
use of bad primes in rational reconstruction}, Math. Comp. 84, 3013--3027 (2015).

\bibitem{BBDLPSS}J.~B\"{o}hm, W.~Decker, S.~Laplagne, G.~Pfister,
A.~Steenpa\ss , S.~Steidel. \emph{Parallel algorithms for normalization}, J.
Symb. Comp. 51, 99--114 (2013).

\bibitem{BDLP1}J.~B\"{o}hm, W.~Decker, G.~Pfister, S.~Laplagne.
\emph{Local to global algorithms for the Gorenstein adjoint ideal of a curve},
Preprint (2015), \href{http://arxiv.org/abs/1505.05040}{arXiv:1505.05040}.

\bibitem{DGPS}
Decker, W., Greuel, G.-M., Pfister, G., Sch{\"o}nemann, H., 2015.
\newblock {\sc Singular} {4-0-2} -- {A} computer algebra system for polynomial
computations. 
\url{http://www.singular.uni-kl.de}



\bibitem{GLS}G.-M.~Greuel, S.~Laplagne, S.~Seelisch, \emph{Normalization
of rings}, J. Symb. Comp. 45(9), 887--901 (2010).

\bibitem{KG}P. Kornerup, R. T. Gregory, \emph{Mapping integers and Hensel
codes onto Farey fractions}, BIT 23, 9--20 (1983).

\bibitem{Mnuk}Mnuk, M.: \emph{An algebraic approach to computing
adjoint curves}. J. Symbolic Comput., 23(2-3), 229-240 (1997).

\bibitem{vanHoeijIntegralBasis}van Hoeij, M.: \emph{An
algorithm for computing an integral basis in an algebraic function field}. J.
Symbolic Comput. 18, no. 4, 353-363 (1994).


\end{thebibliography}
\end{document}